# Two-step benchmarking: Setting more realistically achievable targets in DEA

*Nuria Ramón, José L. Ruiz[*] and Inmaculada Sirvent*

*Centro de Investigación Operativa. Universidad Miguel Hernández. Avd. de la Universidad, s/n 03202-Elche (Alicante), SPAIN*

*E-mail: nramon@umh.es, jlruiz@umh.es, isirvent@umh.es*

*February, 2017*

**Abstract**

The models that set the closest targets have made an important contribution to DEA as tool for the best-practice benchmarking of decision making units (DMUs). These models may help defining plans for improvement that require less effort from the DMUs. However, in practice we often find cases of poor performance, for which closest targets are still unattainable. For those DMUs, we propose a two-step benchmarking approach within the spirit of context-dependent DEA and that of the models that minimize the distance to the DEA efficient frontier. This approach allows to setting more realistically achievable targets in the short term. In addition, it may offer different alternatives for planning improvements directed towards DEA efficient targets, which can be seen as representing improvements in a long term perspective. To illustrate, we examine an example which is concerned with the research performance of public Spanish universities.

## 1. Introduction

In management, organizations use benchmarking for the evaluation of their processes in comparison to best practices of others within a peer group of firms in an industry or sector. In the best-practice benchmarking process the identification of the best firms allows to setting targets, so that organizations may learn from them and develop plans for improvement of some aspects of performance.

Data Envelopment Analysis (DEA) (Charnes et al., 1978) has proven to be a useful tool for the benchmarking of decision making units (DMUs) involved in a production process. In DEA, an empirical production possibility set is constructed from the observations by making some technological assumptions, and the envelopment of such technology determine an efficient frontier formed by the efficient units that is used as reference for the assessments of the remaining DMUs. Nevertheless, Cook et al. (2014) claim that "In the circumstance of benchmarking, the efficient DMUs, as defined by DEA, may not necessarily form a "production frontier", but rather lead to a

[*]Corresponding author.

"best-practice frontier"". Specifically, the points on the best-practice frontier are potential benchmarks for the inefficient units, and the targets are actually the coordinates of these benchmarks, which represent levels of operation for the inefficient DMUs that would make them perform efficiently. As stated in Thanassoulis et al. (2008), in many practical applications one is more interested in determining targets that render the DMUs efficient than in determining their level of inefficiency. See Adler et al. (2013), Zanella et al. (2013) and Dai and Kuosmanen (2014) for some recent references on applications of DEA and benchmarking.

Despite their usefulness for the benchmarking, DEA models often set unrealistic targets, in the sense that they represent plans that are too far from the actual performances. In this respect, the models that minimize the distance from the DMUs to the DEA efficient frontier have made a significant progress over the classical additive DEA models, which maximize the slacks. These models provide a suitable approach to deal with this issue, because they seek to set targets that allow the DMUs to achieve the efficiency with less effort (see Aparicio et al. (2007), Fukuyama et al. (2014a, 2014b), Aparicio and Pastor (2014), Ruiz et al. (2015), Ruiz and Sirvent (2016), Aparicio et al. (2016), Ramón et al. (2016) and Cook et al. (2017)). However, very frequently closest targets are still unachievable for some inefficient DMUs, at least in the short term. This paper focuses on the benchmarking of those units.

Some authors have pointed out these problems of DEA with the benchmarking, and have raised the need of new approaches for planning improvements based on sequential benchmarking, wherein the setting of targets is carried out in several steps. For example, Zhu (2003) states "It is likely that a particular inefficient DMU is unable to immediately improve its performance onto the first level best-practice frontier because of such restrictions as management expertise, available resources, etc. Therefore intermediate (and more easily achievable) targets may be desirable for an inefficient DMU (…) The resulting intermediate targets are *local targets*, whereas the targets on the first level (original) best-practice frontier are *global targets*". Lim et al. (2011) claim "It is very likely to be practically infeasible for an inefficient DMU to achieve the target's efficiency in a single step. That is, if the inefficient DMU is far from the efficient frontier, it will be impossible to reach the frontier in a single move, the more reasonable alternative being to make stepwise gradual improvements in getting to the target". Lozano and Villa (2010) stress that "In Data Envelopment Analysis (DEA) an inefficient unit can be projected onto an efficient target that is far away, i.e., reaching the target may demand large reductions in inputs and increases in outputs. When the inputs and outputs modifications planned are large, it may be troublesome to carry them out all at once." See also the discussion in Kwon et al. (2017) where the authors propose a "better practice" benchmarking approach as distinguished from the traditional "best practice" benchmarking.

A number of approaches have been proposed in DEA for setting targets on a series of nested efficient frontier layers that result from the omission of some DMUs. The idea of using nested efficient frontier layers in DEA is not new. Barr et al. (2000) talk about “peeling the DEA onion” in an approach in which DMUs are separated into a series of nested efficient frontier layers in order to identify a complete ranking of units. Bougnol and Dulá (2006) use that idea in a discussion on the validation of DEA as a ranking tool. Seiford and Zhu (1999) refer to it as “stratification” in the context-dependent DEA approach, which implies the evaluation against an efficient frontier composed by DMUs in a specific performance level (see also Johnson and Zhu (2006)). Within a framework of benchmarking, Lim et al. (2011) establish a sequence of targets on the nested efficient frontier layers that result from a benchmark selection according to criteria of attractiveness, progress and infeasibility as introduced in the context-dependent DEA approach in Seiford and Zhu (2003); Kwon et al. (2017) also use this “stratification” in their “better practice” benchmarking approach, which combines DEA (radial models) and neural networks. Other sequential DEA approaches for benchmarking include Lozano and Villa (2005, 2010), which propose strategies of gradual improvements with successive, intermediate targets; Fang (2015), which uses a similar idea for centralized models; and Estrada et al. (2009), which propose a stepwise benchmarking based on input similarities by using a method that combines DEA and self-organizing maps (SOM). See also Petrovic et al. (2014) for a multi-criteria decision making (MCDM) approach based on the multi-level outranking ELECTRE method.

In the present paper, we propose a two-step benchmarking determined by a couple of nested DEA efficient frontiers, following an approach based on closest targets. Thus, the distinctive feature with respect to other stepwise benchmarking approaches based on context-dependent DEA is that we use models that minimize the distance to the efficient frontiers associated with the different levels of performance. That is, special attention is paid to the fact that improvements should require as little effort as possible from the DMUs. In addition, we also seek to find targets that are efficient in the Pareto sense. This is why our approach is developed within the framework of non-radial models. As Thanassoulis et al. (2008) state, non-radial models are the appropriate instrument for the setting of targets and the benchmarking.

Specifically, the purpose of the approach proposed here is twofold: 1) setting more realistically achievable targets in the short term for the inefficient DMUs, and 2) establishing a sequential plan of improvements for the inefficient DMUs directed towards meeting long terms targets. In order to do so, we consider the following two DEA efficient frontiers: the frontier associated with the whole set of DMUs and the one resulting from the omission of the DMUs initially rated as efficient. As in context-dependent DEA, this “stratification” allows to classifying the DMUs into different levels of performance. The DMUs in the outer frontier are in a first level of performance, those in the other frontier are in a second level of performance, and the remaining DMUs are in a third level of performance (or higher, if the "stratification" process is continued). The DMUs in the second level of performance can be benchmarked against the 1st-level efficient frontier, while the remaining inefficient DMUs can be benchmarked against the 1st- or 2nd-level efficient frontier. However, for the inefficient DMUs in the 3rd level (or higher), the targets set on the original DEA efficient frontier are often unattainable, even though they result from a model that minimizes the distance to such frontier. This is why targets on that frontier for those units can be considered as representing improvements in a long term perspective (note that in DEA the DMUs are assumed to be homogeneous, that is, that they all are comparable in terms of the variables chosen). We design a benchmarking process directed towards long term targets which, in a first step, allows to setting more realistically achievable targets on the 2nd-level efficient frontier, which can be seen as a best-practice frontier determined by DMUs in a more similar level of performance to the units under evaluation than those that form the original DEA efficient frontier. In addition, this sequential approach allows also to distributing efforts in the way for improvement, which is especially useful for inefficient DMUs that show a poor performance. It should be pointed out that the ultimate targets set by this approach on the outer frontier might be further from the DMU under evaluation than its closest targets on that frontier, but the possible extra effort that could be needed to reach them would be justified by the fact that improvements would be made in two steps.

The benchmarking model proposed consider two objectives for setting the sequence of targets corresponding to a given DMU: minimizing the gap between actual performances and intermediate targets and minimizing the gap between these latter and the long term targets. These two objectives are combined by means of a parameter which allows us to adjust the importance attached to each of them. Through the specification of this parameter, several sequential strategies for improvement can be defined, in the sense that each sequence of targets establishes the way in which efforts are distributed in the two steps regarding both the specific inputs/outputs to act upon and the magnitude of the changes to be accomplished. Thus, each DMU can establish the plan for improvements to implement making a choice among those alternatives according its circumstances.

The paper unfolds as follows. In section 2 we develop the two-step benchmarking model that sets a sequence of targets for inefficient DMUs. The approach proposed is illustrated in section 3 with an example in which research performance of public Spanish universities is examined. The final section presents conclusions.

## 2. The two-step benchmarking model

Throughout the paper we suppose that we have *n* DMUs which use *m* inputs to produce *s* outputs. These are denoted by $\left(X_j, Y_j\right)$, $j = 1,\ldots,n$. It is assumed that $X_j = \left(x_{1j},\ldots,x_{mj}\right) \geq 0$, $X_j \neq 0$, $j = 1,\ldots,n$, and $Y_j = \left(y_{1j},\ldots,y_{sj}\right) \geq 0$, $Y_j \neq 0$, $j = 1,\ldots,n$. For purposes of benchmarking, it is also assumed that the production possibility set (PPS) $T = \left\{\left(X,Y\right) / \text{ X can produce Y}\right\}$ satisfies the classical postulates: convexity, variable returns to scale (VRS) and free disposability. Therefore, T can be characterized as follows $T = \left\{\left(X,Y\right) \Big/ X \geq \sum_{j=1}^{n} \lambda_j X_j,\ Y \leq \sum_{j=1}^{n} \lambda_j Y_j,\ \sum_{j=1}^{n} \lambda_j = 1,\ \lambda_j \geq 0\right\}$ (see Banker et al., 1984).

For the development of the two-step benchmarking we consider two DEA efficient frontiers. Let $\partial(T)$ be the Pareto-efficient frontier formed by the efficient DMUs of T. We denote by E the set consisting of those DMUs, which are in a first level of performance. Likewise, if we denote by $T^I$ the set of DMUs in T resulting from the omission of those in E, then we can consider the 2nd-level efficient frontier $\partial(T^I)$ (in the sense of Pareto), which is the frontier formed by the efficient DMUs of $T^I$ (denote by $E^I$ that set of efficient units). Obviously, both frontiers are nested.

The DMUs in a second level of performance can be benchmarked against those in the first level. Therefore, for each $DMU_0$ in $\partial(T^I)$, targets $\left(\hat{X}_0^*, \hat{Y}_0^*\right)$ can be set by using a model that minimizes the distance to $\partial(T)$ like the following

$$
\begin{aligned}
\text{Min} \quad & \left\| \left(X_0, Y_0\right) - \left(\hat{X}_0, \hat{Y}_0\right) \right\|_1^{\omega} \\
\text{s.t.:} \quad & \\
& \left(\hat{X}_0, \hat{Y}_0\right) \in \partial(T) \\
& \left(\hat{X}_0, \hat{Y}_0\right) \text{ dominates } \left(X_0, Y_0\right)
\end{aligned}
\tag{1}
$$

wherein the weighted $L_1$-norm $\left\| \left(X_0, Y_0\right) - \left(\hat{X}_0, \hat{Y}_0\right) \right\|_1^{\omega} = \sum_{i=1}^{m} \frac{x_{i0} - \hat{x}_{i0}}{x_{i0}} + \sum_{r=1}^{s} \frac{\hat{y}_{r0} - y_{r0}}{y_{r0}}$ is used. Model (1) set the closest targets to $DMU_0$ on the 1st-level efficient frontier.

Taking into account the characterization of $\partial(T)$ in terms of a set of linear constraints (see, for instance, Aparicio et al. (2007)), we have the following operative formulation of model (1), which is expressed in terms of the usual slacks

$$
\begin{aligned}
\text{Min} \quad & \sum_{i=1}^{m} \frac{s_{i0}^{E-}}{x_{i0}} + \sum_{r=1}^{s} \frac{s_{r0}^{E+}}{y_{r0}} \\
\text{s.t.:} \quad & \\
& \sum_{j \in E} \lambda_j^E x_{ij} = x_{i0} - s_{i0}^{E-} && i = 1,\dots,m \\
& \sum_{j \in E} \lambda_j^E y_{rj} = y_{r0} + s_{r0}^{E+} && r = 1,\dots,s \\
& \sum_{j \in E} \lambda_j^E = 1 \\
& -\sum_{i=1}^{m} v_i^E x_{ij} + \sum_{r=1}^{s} u_r^E y_{rj} + u_0^E + d_j^E = 0 && j \in E \\
& v_i^E x_{i0} \geq 1 && i = 1,\dots,m \\
& u_r^E y_{r0} \geq 1 && r = 1,\dots,s \\
& d_j^E \leq M^E b_j^E && j \in E \\
& \lambda_j^E \leq M^E (1 - b_j^E) && j \in E \\
& s_{i0}^{E-} \geq 0 && i = 1,\dots,m \\
& s_{r0}^{E+} \geq 0 && r = 1,\dots,s \\
& \lambda_j^E, d_j^E \geq 0 && j \in E \\
& b_j^E \in \{0,1\} && j \in E \\
& u_0^E \quad \text{free}
\end{aligned}
\tag{2}
$$

$M^E$ being a big positive quantity. Therefore, by using the optimal solutions of (2), targets $\left(\hat{X}_0^*, \hat{Y}_0^*\right)$ can be found as follows

$$
\begin{aligned}
\hat{X}_0^* &= X_0 - S_0^{E-*} \left( = \sum_{j \in E} \lambda_j^{E*} X_j \right) \\
\hat{Y}_0^* &= Y_0 + S_0^{E+*} \left( = \sum_{j \in E} \lambda_j^{E*} Y_j \right)
\end{aligned}
\tag{3}
$$

where $S_0^{E-*} = \left(s_{10}^{E-*}, \dots, s_{m0}^{E-*}\right)'$ and $S_0^{E+*} = \left(s_{10}^{E+*}, \dots, s_{s0}^{E+*}\right)'$.

For each of the remaining inefficient units, say $DMU_0$, the approach proposed seeks to find a sequence of targets $\left(\hat{X}_0^{I*}, \hat{Y}_0^{I*}\right)$, on the intermediate frontier, and $\left(\hat{X}_0^{*}, \hat{Y}_0^{*}\right)$, on the outer frontier, which allow to establishing a sequential plan of improvements directed towards best practices as determined by DMUs in a first level of performance. Such plan is expected to be implemented through more realistic efforts, especially in the case of units that show a poor performance. In addition to lying on their respective efficient frontiers, in the determination of such targets we seek both minimizing the gap between $\left(X_0, Y_0\right)$ and $\left(\hat{X}_0^{I*}, \hat{Y}_0^{I*}\right)$ and minimizing the gap between $\left(\hat{X}_0^{I*}, \hat{Y}_0^{I*}\right)$ and $\left(\hat{X}_0^{*}, \hat{Y}_0^{*}\right)$. Again, the idea is that targets define directions of improvement that requires as little effort as possible from $DMU_0$. In terms of the weighted $L_1$-norm, we should consider minimizing both

$$\begin{aligned} &\left\|\left(X_0, Y_0\right)-\left(\hat{X}_0^{I}, \hat{Y}_0^{I}\right)\right\|_1^{\omega}=\sum_{i=1}^{m}\frac{x_{i0}-\hat{x}_{i0}^{I}}{x_{i0}}+\sum_{r=1}^{s}\frac{\hat{y}_{r0}^{I}-y_{r0}}{y_{r0}} \qquad (4.1)\\ &\left\|\left(\hat{X}_0^{I}, \hat{Y}_0^{I}\right)-\left(\hat{X}_0, \hat{Y}_0\right)\right\|_1^{\omega}=\sum_{i=1}^{m}\frac{\hat{x}_{i0}^{I}-\hat{x}_{i0}}{x_{i0}}+\sum_{r=1}^{s}\frac{\hat{y}_{r0}-\hat{y}_{r0}^{I}}{y_{r0}} \qquad (4.2)\end{aligned} \tag{4}$$

where all the deviations are assumed to be non-negative in order to define a strategy of continuous improvements.

When the benchmarking is carried out against $\partial(T)$ in a single step, considering other referents for $DMU_0$ different from its closest targets might have no sense (there should be some good reason for making an extra effort). However, focusing only on its closest targets on $\partial(T^I)$ in the first step of a sequential benchmarking might not be an appropriate strategy, because these targets could be very far from $\partial(T)$, thus leading to a second step of the plan for improvements that would require much effort. For this reason, in the two-step benchmarking other points on $\partial(T^I)$ aside from the closest targets should be considered.

Ideally, one would like to minimize simultaneously both (4.1) and (4.2). However, it is obvious that the targets $\left(\hat{X}_0^{I}, \hat{Y}_0^{I}\right)$ that minimize (4.1) will not be necessarily those that minimize (4.2), and the other way around too. As a compromise between these two objectives, we propose to minimize a convex combination of them

$$\alpha\left\|\left(X_0, Y_0\right)-\left(\hat{X}_0^{I}, \hat{Y}_0^{I}\right)\right\|_1^{\omega}+\left(1-\alpha\right)\left\|\left(\hat{X}_0^{I}, \hat{Y}_0^{I}\right)-\left(\hat{X}_0, \hat{Y}_0\right)\right\|_1^{\omega}, \tag{5}$$

where $0 \le \alpha \le 1$. This will provide us with non-dominated solutions. Through the specification of $\alpha$, we may adjust the importance to be attached to each of the two objectives. For example, setting $\alpha = 1$ means to look for the closest intermediate targets, while with $\alpha = 0$ the search is directed towards intermediate targets that are as close as possible to the 1st-level efficient frontier. Values of $\alpha$ in (0,1) will lead to plans in between those two extremes. A series of sequential targets can then be generated for different values of $\alpha$. Thus, different performance improvement strategies oriented to achieve long term targets can be suggested to $DMU_0$, so that this unit may choose among them the best plan to implement according to its circumstances.

Bearing in mind the above, a model that provides the sequence of targets that is sought for $DMU_0$ (for a given $\alpha$) can be formulated as follows

$$
\begin{aligned}
\text{Min} \quad & \alpha \left\| \left(X_0, Y_0\right) - \left(\hat{X}_0^I, \hat{Y}_0^I\right) \right\|_1^{\omega} + (1-\alpha) \left\| \left(\hat{X}_0^I, \hat{Y}_0^I\right) - \left(\hat{X}_0, \hat{Y}_0\right) \right\|_1^{\omega} \\
\text{s.t.:} \quad & \\
& \left(\hat{X}_0^I, \hat{Y}_0^I\right) \in \partial(T^I) \\
& \left(\hat{X}_0, \hat{Y}_0\right) \in \partial(T) \\
& \left(\hat{X}_0^I, \hat{Y}_0^I\right) \text{ dominates } \left(X_0, Y_0\right) \\
& \left(\hat{X}_0, \hat{Y}_0\right) \text{ dominates } \left(\hat{X}_0^I, \hat{Y}_0^I\right)
\end{aligned}
\tag{6}
$$

Using again the results in Aparicio et al. (2007), we may derive an operative formulation of model (6) as a mixed-integer linear programming problem

$$
\begin{aligned}
&\text{Min} \;\; \alpha\left(\sum_{i=1}^{m}\frac{s_{i0}^{I-}}{x_{i0}}+\sum_{r=1}^{s}\frac{s_{r0}^{I+}}{y_{r0}}\right)+(1-\alpha)\left(\sum_{i=1}^{m}\frac{s_{i0}^{E-}}{x_{i0}}+\sum_{r=1}^{s}\frac{s_{r0}^{E+}}{y_{r0}}\right)\\
&\text{s.t.:}\\
&\sum_{j\in E^I}\lambda_j^I x_{ij}=x_{i0}-s_{i0}^{I-} && i=1,\ldots,m\\
&\sum_{j\in E^I}\lambda_j^I y_{rj}=y_{r0}+s_{r0}^{I+} && r=1,\ldots,s\\
&\sum_{j\in E^I}\lambda_j^I=1\\
&-\sum_{i=1}^{m}v_i^I x_{ij}+\sum_{r=1}^{s}u_r^I y_{rj}+u_0^I+d_j^I=0 && j\in E^I\\
&v_i^I x_{i0}\geq 1 && i=1,\ldots,m\\
&u_r^I y_{r0}\geq 1 && r=1,\ldots,s\\
&d_j^I\leq M^I b_j^I && j\in E^I\\
&\lambda_j^I\leq M^I(1-b_j^I) && j\in E^I\\
&\sum_{j\in E}\lambda_j^E x_{ij}=\left(x_{i0}-s_{i0}^{I-}\right)-s_{i0}^{E-} && i=1,\ldots,m\\
&\sum_{j\in E}\lambda_j^E y_{rj}=\left(y_{r0}+s_{r0}^{I+}\right)+s_{r0}^{E+} && r=1,\ldots,s\\
&\sum_{j\in E}\lambda_j^E=1\\
&-\sum_{i=1}^{m}v_i^E x_{ij}+\sum_{r=1}^{s}u_r^E y_{rj}+u_0^E+d_j^E=0 && j\in E\\
&v_i^E x_{i0}\geq 1 && i=1,\ldots,m\\
&u_r^E y_{r0}\geq 1 && r=1,\ldots,s\\
&d_j^E\leq M^E b_j^E && j\in E\\
&\lambda_j^E\leq M^E(1-b_j^E) && j\in E\\
&s_{i0}^{I-},s_{i0}^{E-}\geq 0 && i=1,\ldots,m\\
&s_{r0}^{I+},s_{r0}^{E+}\geq 0 && r=1,\ldots,s\\
&b_j^I\in\{0,1\},\, j\in E^I,\, b_j^E\in\{0,1\},\, j\in E\\
&u_0^I,u_0^E \quad \text{free}\\
&\lambda_j^I,d_j^I\geq 0,\, j\in E^I,\,\lambda_j^E,d_j^E\geq 0,\, j\in E && (7)
\end{aligned}
$$

where $M^E$ and $M^I$ are two big positive quantities. Therefore, for every $\alpha$, the optimal solutions of (7) provide the sequence of targets $\left(\hat{X}_0^{I*},\hat{Y}_0^{I*}\right)$

$$\begin{aligned} \hat{X}_0^{I*} &= X_0 - S_0^{I-*}\left(=\sum_{j\in E^I}\lambda_j^{I*}X_j\right) \\ \hat{Y}_0^{I*} &= Y_0 + S_0^{I+*}\left(=\sum_{j\in E^I}\lambda_j^{I*}Y_j\right) \end{aligned} \tag{8}$$

as intermediate targets, where $S_0^{I-*}=\left(s_{10}^{I-*},\ldots,s_{m0}^{I-*}\right)'$ and $S_0^{I+*}=\left(s_{10}^{I+*},\ldots,s_{s0}^{I+*}\right)'$, and $\left(\hat{X}_0^{*},\hat{Y}_0^{*}\right)$

$$\begin{aligned} \hat{X}_0^{*} &= \left(X_0 - S_0^{I-*}\right) - S_0^{E-*}\left(=\sum_{j\in E}\lambda_j^{E*}X_j\right) \\ \hat{Y}_0^{*} &= \left(Y_0 + S_0^{I+*}\right) + S_0^{E+*}\left(=\sum_{j\in E}\lambda_j^{E*}Y_j\right) \end{aligned} \tag{9}$$

as targets on the outer frontier, where $S_0^{E-*}=\left(s_{10}^{E-*},\ldots,s_{m0}^{E-*}\right)'$ and $S_0^{E+*}=\left(s_{10}^{E+*},\ldots,s_{s0}^{E+*}\right)'$.

Remark 1. We note that, for $\alpha=1$, the targets on the 1st-level efficient frontier provided by (7) are found without following any criterion of optimality. For this reason, in that case we propose to set them as the optimal solution of (2) when $\left(\hat{X}_0^{I*},\hat{Y}_0^{I*}\right)$, as obtained from (7) for $\alpha=1$, is evaluated. That is, setting $\left(\hat{X}_0^{*},\hat{Y}_0^{*}\right)$ as the closest targets to $\left(\hat{X}_0^{I*},\hat{Y}_0^{I*}\right)$ on the outer frontier.

Remark 2. Constraints $d_j^I \le M^I b_j^I,\ j\in E^I$, and $\lambda_j^I \le M^I(1-b_j^I),\ j\in E^I$, which include the classical big M and binary variables, seek that the DMU$_j$'s in $E^I$ that participate actively as a referent in the evaluation of DMU$_0$ necessarily belong all to a supporting hyperplane containing a facet of $\partial(T^I)$. Nevertheless, (7) can be solved in practice by reformulating these constraints using Special Ordered Sets (SOS) (Beale and Tomlin, 1970), which avoid the need to specify $M^I$. SOS Type 1 is a set of variables where at most one variable may be nonzero. Therefore, if we remove these two groups of constraints from the formulation and define instead a SOS Type 1 for each pair of variables $\left\{\lambda_j^I, d_j^I\right\}$, $j\in E^I$, then it is ensured that $\lambda_j^I$ and $d_j^I$ cannot be simultaneously positive for DMU$_j$'s, $j\in E^I$. CPLEX Optimizer (and also LINGO) can solve LP problems with SOS. SOS variables have already used for solving models like (7) in Ruiz and Sirvent (2016), Aparicio et al. (2016) and Cook et al.

(2017). This reasoning also applies to the constraints $d_j^E \leq M^E b_j^E$, $j \in E$, and $\lambda_j^E \leq M^E(1-b_j^E)$, $j \in E$, in the same model.

Remark 3. The proposed approach can be straightforwardly extended to consider more than two nested DEA frontier, if the benchmarking against DMUs in more than two levels of performance could be of help for some inefficient unit. Nevertheless, carrying out a plan for improvements that needs to be implemented in many steps (before the next period of performance) may appear unrealistic in practice.

Remark 4. It should be highlighted that this two-step benchmarking may lead to points on $\partial(T)$ which are not the closest targets to $DMU_0$ on that frontier, in return for allowing the implementation of a plan for improvements in several steps. Note, in any case, that we cannot force $DMU_0$ to go to its closest targets on $\partial(T)$ following a sequential approach like the one proposed here, because it cannot be ensured that both points can always be connected through $\partial(T^I)$ following a plan of continuous improvements (this will be shown graphically in the numerical example).

Numerical example

Consider the DMUs in Table 1 which produce two outputs by using one constant input. A conventional DEA analysis reveals that A, B, C and D are the extreme efficient DMUs. Thus, they form the 1st-level efficient frontier, $\partial(T)$. Among the remaining units, we highlight the cases of DMUs 4 and 5, which show a poor performance. They are actually in a 3rd (or higher) level of performance. The approach proposed here is of particular interest for the DMUs in that situation. A model that minimizes distance to the frontier would project DMUs 4 and 5 on to the points 4' and 5' on $\partial(T)$ (see Figure 1). Despite such models set the closest targets, in these particular cases they would require much effort from those units to perform efficiently. Specifically, DMUs 4 and 5 would have to produce 8.25 and 8.625 units of $y_2$, when their actual production is 4 and 5 units, respectively, while maintaining the actual level of performance in $y_1$. This situation could be suggesting the need of a little more realistic plan for improvements.

Table 1

The two-step benchmarking we propose can be of help for setting more realistically achievable targets in the short term. It may also help to define a sequential strategy directed towards targets on $\partial(\mathrm{T})$. As has been explained, that approach considers the 2nd-level efficient frontier, which in this case is the one associated with the set of units resulting from the omission of DMUs A, B, C and D. Such frontier, $\partial(\mathrm{T}^{\mathrm{I}})$, is determined by DMUs 1, 2 and 3 (see again Figure 1). Model (7) may offer different alternatives for planning improvements through the specification of $\alpha$. Table 2 records the series of sequential targets for values of $\alpha = 0,\ 0.25, 0.5, 0.75, 1$.

Table 2

Consider the case of DMU 4. For $\alpha = 0.25,\ 0.5$ and $0.75$, the optimal solutions of model (7) show how to distribute the effort to achieve the closest targets on $\partial(\mathrm{T})$ (previously mentioned) through the 2nd-level best-practice frontier $\partial(\mathrm{T}^{\mathrm{I}})$. Specifically, by increasing $y_2$ from 4 to 6.75 in a first step and leaving for a second stage the remaining effort up to 8.25 units. However, other $\alpha$'s give rise to other improvement alternatives, which, in spite of possibly requiring quantitatively a little more effort, this DMU might consider. For example, when $\alpha = 1$, model (7) proposes a plan that, unlike the one corresponding to the specifications $\alpha = 0.25,\ 0.5$ and $0.75$ just discussed, involves improvements only in $y_1$. Specifically, it suggest to increase this variable from 5 to 10 in two steps by making the same effort in each of them (to eventually projecting on 4''). In addition, for $\alpha = 0$, we find a new plan which proposes to improve simultaneously both outputs in the two steps using DMU 2 as the benchmark in the first step. Thus, this sequential benchmarking approach offers different alternatives for improvement, so that DMU 4 can establish the plan to implement by making a choice among them taking into account its individual circumstances.

Something similar happens to DMU 5. Note however that for that unit it would be impossible to implement a sequential planning ending in the closest targets on $\partial(\mathrm{T})$ (point 5') following an approach of continuous improvements (in which the intermediate benchmark is Pareto efficient), because there is no path from DMU 5 to 5' through $\partial(\mathrm{T}^{\mathrm{I}})$ satisfying that requirement.

Figure 1

## 3. Empirical illustration

As an illustration of the proposed approach, in this section we apply it to the evaluation of research activities of public Spanish universities. Research is one of main areas of performance of universities, because of its importance for the social development and the improvement of the quality of life of people. The public Spanish universities may be seen as a set of homogeneous DMUs that undertake similar activities and produce comparable results regarding research, so that a common set of outputs can be defined for their analysis. Specifically, we consider the following inputs and outputs, which relate as usual to human and physical capital on one hand, and publications (regarding both production and quality) and incomes on the other:

INPUTS

- ASTAFF: Full-time equivalent academic staff (academic year 2013-14).
- EXPEND: Expenditure exactly accounts for staff expenditure, expenditure on goods and services, financial expenditure and current transfers. EXPEND is an indicator traditionally used for comparing institutions, because it reflects the budgetary effort made by the universities in the delivery of their activities (year 2013).

OUTPUTS

- ARTICLES: Total production of articles published in the period 2009-13.
- %Q1: Percentage of articles published in journals classified in the 1st quartile of the corresponding category.
- INCOMES: Total income for research and development activities (year 2013).
- THESIS: Number of thesis defended in the period 2009-13.

Data for these variables have been taken from the corresponding reports by the Conference of Rectors of the Spanish Universities (CRUE) and by the Foundation Knowledge and Development (FCyD). Our study here can thus complement the results provided by those reports with an evaluation of the universities from a perspective focused on benchmarking and target setting.

The sample consists of 46 (out of 48) public Spanish universities. Table 3 reports a descriptive summary. An initial DEA analysis reveals 13 universities as technically efficient.

Table 3

Table 4 records the results derived from model (2) regarding the setting of targets for some inefficient universities that show a poor performance, as representative cases: University of Cadiz (UCA), University of Las Palmas de Gran Canaria (ULPGC), University of Salamanca (USAL) and University Carlos III (UCAR). We can see that, although this model minimizes the distance to the efficient frontier, the universities would have to make a strong effort if they would want to achieve those targets in a single step. For example, UCA would have to improve simultaneously its use of resources (both in staff and expenditures) and its publications (both in quantity and in quality) by more than 30% in most cases. Something similar occurs with ULPGC and USAL, while UCAR should make a strong effort to achieve its targets in ASTAFF, EXPEND, %Q1 and THESIS.

Table 4

Table 5 reports the results of the sequential benchmarking provided by (7) for those universities for $\alpha = 0,\ 0.25,\ 0.5,\ 0.75 \text{ and } 1$. Such analysis is carried out by using the 2nd-level efficient frontier formed by the efficient universities in the sample resulting from eliminating the 13 efficient universities in the 1st level of performance. Universities UCA, ULPGC, USAL and UCAR are inefficient even with respect to such intermediate frontier, that is, they are in a 3rd (or higher) level of performance. Let us consider the case of UCA. For $\alpha = 0.5,\ 0.75 \text{ and } 1$ this approach sets intermediate targets that allow this university to distribute the effort to meet the targets provided by (2) (those in Table 4). Specifically, the targets in rows $\alpha = 0.75 \text{ and } 1$ suggest a strategy based on making first an important effort in the use of resources (both in ASTAFF and EXPEND) and leaving for a later stage the bulk of actions oriented to improve publications (both ARTICLES and %Q1). Something similar occurs with UCAR. In that case, the results provided by model (7) suggest to make the improvements regarding THESIS in a first stage and those corresponding to the resources in a second step, while at the same time it is suggested to improve the quality of publications (%Q1) gradually.

For ULPGC, rows $\alpha = 0.5 \text{ and } 0.75$ show again how to distribute the efforts needed to achieve the targets set by (2). Basically, improving first the use of resources and then the production of ARTICLES. In addition, this approach offers an alternative strategy when $\alpha = 1$, which is globally more demanding, but which ULPGC might consider. Such plan is more demanding in ASTAFF and less in ARTICLES, in exchange for suggesting a little more effort in %Q1.

Finally, for USAL, we note that none of the plans for improvements associated with the different $\alpha$'s end in the closest targets reported in Table 4. We highlight the plan associated with $\alpha = 0.75$, because it is the one which distributes more the efforts between the two steps. Nevertheless, more effort should be made in the first stage, leaving for a second step some improvements in EXPEND and ARTICLES.

Table 5

## 4. Conclusions

DEA is a benchmarking tool that assumes the homogeneity of the DMUs involved in the evaluation of their performance. Therefore, inefficient DMUs should, in principle, agree with the targets set by the DEA models for them (especially when they set the closest targets), because the DMUs are all comparable. However, some inefficient DMUs show frequently a very poor performance, so DEA targets in those cases become unattainable, at least in the short term. For the benchmarking of those DMUs, we propose a sequential approach directed towards DEA targets, which can be seen as representing improvements in a longer term perspective. At the same time, it allows, in a first step, to setting more realistically achievable targets on an intermediate best-practice frontier formed by DMUs in a more similar level of performance. This approach has been developed within the framework of models that minimize the distance to the Pareto-efficient frontier, because we seek to show the DMUs the directions for improvement that require less effort. The example examined has also shown its usefulness in practice in the sense that it may offer different alternatives on how to distribute efforts to be made for improvement. Thus, managers of organizations can establish their planning making a choice among them taking into account their individual circumstances.

## Acknowledgments

This research has been supported through Grant MTM2016-76530-R (AEI/FEDER, UE).

## References

Adler, N., Liebert, V. & Yazhemsky, E. (2013). "Benchmarking airports from a managerial perspective", *Omega*, 41(2), 442-458.

Aparicio, J. & Pastor, J.T. (2014). "Closest targets and strong monotonicity on the strongly efficient frontier in DEA", *Omega*, 44, 51–57.

Aparicio, J., Cordero, J.M & Pastor, J.T. (2016). "The determination of the least distance to the strongly efficient frontier in Data Envelopment Analysis oriented models: Modelling and computational aspects", *Omega* (2016), http://dx.doi.org/10.1016/j.omega.2016.09.008.

Aparicio, J., Ruiz, J.L. & Sirvent, I. (2007). "Closest targets and minimum distance to the Pareto-efficient frontier in DEA", *Journal of Productivity Analysis*, 28(3), 209–218, 2007.

Banker, R.D., Charnes, A & Cooper, W.W. (1984). "Some models for estimating technical and scale inefficiencies in data envelopment analysis", *Management Science*, 30(9), 1078–1092

Barr, R.S., Durchholz, M.L. & Seiford, L. (2000). "Peeling the DEA onion: Layering and rank-ordering DMUs using tiered DEA", *Southern Methodist University Technical Report*, 5, 1-24.

Beale, E.M.L. & Tomlin, J.A. (1970). "Special facilities in a general mathematical programming system for non-convex problems using ordered sets of variables", in J. Lawrence (ed.), *Proceedings of the Fifth International Conference on Operational Research* (Tavistock Publications, London, 1970), 447-454.

Bougnol, M.L. & Dulá, J.H. (2006). "Validating DEA as a ranking tool: An application of DEA to assess performance in higher education", *Annals of Operations Research*, 145, 339-365.

Charnes, A., Cooper, W.W. & Rhodes, E. (1978). "Measuring the efficiency of decision making units", *European Journal of Operational Research*, 2(6), 429-444.

Cook W.D, Tone K & Zhu, J. (2014). "Data envelopment analysis: prior to choosing a model", *Omega*, 44(C), 1–4

Cook, W.D., Ruiz, J.L., Sirvent, I. & Zhu, J. (2017). "Within-group common benchmarking using DEA", *European Journal of Operational Research*, 256, 901-910.

Dai, X. & Kuosmanen, T. (2014). "Best-practice benchmarking using clustering methods: Application to energy regulation", *Omega*, 42(1), 179-188.

Estrada, S.A., Song, H.S., Kim. Y.A., Namn, S.H. & Kang, S.C. (2009). "A method of stepwise benchmarking for inefficient DMUs based on the proximity-based target selection", *Expert Systems with Applications*, 36(9), 11595-11604.

Fang, L. (2015). "Centralized resource allocation based on efficiency analysis for step-by-step improvement paths", *Omega*, 51, 24-28.

Fukuyama, H., Maeda, Y., Sekitani, K & Shi, J. (2014a). "Input-output substitutability and strongly monotonic p-norm least-distance DEA measures", *European Journal of Operational Research,* 237(3), 997–1007.

Fukuyama, H., Maeda, Y., Sekitani, K & Shi, J. (2014b). "Distance optimization approach to ratio-form efficiency measures in data envelopment analysis", *Journal of Productivity Analysis,* 42, 175–186.

Kwon, H.B., Marvel, J.H. & Roh, J.J. (2017). "Three-stage performance modeling using DEA-BPNN for better practice benchmarking", *Expert Systems with Applications*, 71, 429-441.

Johnson, S.A. & Zhu, J. (2006). "Identifying ''best'' applicants in recruiting using data envelopment analysis", *Socio-Economic Planning Sciences*, 37, 125-139.

Lim, S., Bar, H. & Lee, L.H. (2011). "A study on the selection of benchmarking paths in DEA", *Expert Systems with Applications,* 38(6), 7665-7673.

Lozano, S. & Villa, G. (2005). "Determining a sequence of targets in DEA", *Journal of the Operational Research Society,* 56(12), 1439-1447.

Lozano, S. & Villa, G. (2010). "Gradual technical and scale efficiency improvement in DEA", *Annals of Operations Research*, 56(12), 1439-1447.

Petrovic, M., Bojkovic, N., Anic, I., Stamenkovic, M. & Tarle, S.P. (2014). "An electre-based decision aid tool for stepwise benchmarking: An application over EU Digital Agenda targets", *Decision Support Systems*, 59, 230-241.

Ramón, N., Ruiz, J.L. & Sirvent, I. (2016). "On the Use of DEA Models with Weight Restrictions for Benchmarking and Target Setting". In Aparicio, J., Lovell, C.A.K. and Pastor, J.T. (eds.), *Advances in Efficiency and Productivity,* Springer.

Ruiz, J.L. & Sirvent, I. (2016). "Common benchmarking and ranking of units with DEA", *Omega,* 65, 1-9.

Ruiz, J.L., Segura, J.V., & Sirvent, I. (2015). "Benchmarking and target setting with expert preferences: An application to the evaluation of educational performance of Spanish universities", *European Journal of Operational Research,* 242(2), 594-605.

Seiford, L.M. & Zhu, J. (1999). "Profitability and Marketability of the Top 55 U.S. Commercial Banks", *Management Science*, 45(9), 1270-1288.

Seiford, L.M. & Zhu, J. (2003). "Context-dependent data envelopment analysis—Measuring attractiveness and progress", *Omega,* 31, 397-408.

Thanassoulis, E., Portela, M.C.A.S., & Despić, O. (2008). "Data Envelopment Analysis: The mathematical programming approach to efficiency analysis". In Harold, O. Fried, C.A. Knox

Lovell & Shelton S. Schmidt (Eds.), *The measurement of productive efficiency and productivity growth*, New York, NY (US): Oxford University Press.

Zanella, A., Camanho, A.S. & Dias, T.G. (2013). "Benchmarking countries' environmental performance", *Journal of the Operational Research Society*, 64(3), 426-438.

Zhu, J. (2003). "Quantitative models for performance evaluation and benchmarking: data envelopment analysis with spreadsheets", Springer International Series.

| DMU | $x_1$ | $y_1$ | $y_2$ |
|---|---|---|---|
| A | 1 | 2 | 9 |
| B | 1 | 6 | 8 |
| C | 1 | 9 | 6 |
| D | 1 | 10.5 | 3 |
| 1 | 1 | 3.5 | 1.5 |
| 2 | 1 | 6.5 | 6.5 |
| 3 | 1 | 8.5 | 1.5 |
| 4 | 1 | 5 | 4 |
| 5 | 1 | 3 | 5 |

Table 1. Data of numerical example

| | 1st step | | | 2nd step | | |
|---|---|---|---|---|---|---|
| | x | $y_1$ | $y_2$ | x | $y_1$ | $y_2$ |
| DMU 4 | 1 | 5 | 4 | | | |
| α=1 | 1 | 7.5 | 4 | 1 | 10 | 4 |
| α=0.75 | 1 | 5 | 6.75 | 1 | 5 | 8.25 |
| α=0.5 | 1 | 5 | 6.75 | 1 | 5 | 8.25 |
| α=0.25 | 1 | 5 | 6.75 | 1 | 5 | 8.25 |
| α=0 | 1 | 6.5 | 6.5 | 1 | 6.5 | 7.667 |
| DMU 5 | 1 | 3 | 5 | | | |
| α=1 | 1 | 3.5 | 7 | 1 | 3.5 | 8.625 |
| α=0.75 | 1 | 3.5 | 7 | 1 | 3.5 | 8.625 |
| α=0.5 | 1 | 3.5 | 7 | 1 | 3.5 | 8.625 |
| α=0.25 | 1 | 3.5 | 7 | 1 | 3.5 | 8.625 |
| α=0 | 1 | 6.5 | 6.5 | 1 | 6.5 | 7.667 |

Table 2. Results of numerical example

| | ASTAFF | EXPEND | ARTICLES | %Q1 | INCOMES | THESIS |
|---|---|---|---|---|---|---|
| Mean | 1706.6 | 15052.79 | 6652.8 | 51.40 | 2105.21 | 962.8 |
| Standard Dev. | 1043.4 | 9750.43 | 4965.5 | 5.91 | 1773.82 | 838.9 |
| Minimum | 355.1 | 3643.33 | 1036 | 39.45 | 209.43 | 159 |
| Maximum | 4864.6 | 44645.32 | 21441 | 62.29 | 7299.48 | 3959 |

Table 3. Descriptive summary

| University | | ASTAFF | EXPEND | ARTICLES | %Q1 | INCOMES | THESIS |
|---|---|---|---|---|---|---|---|
| UCA | Data | 1392.5 | 10923.17 | 2546 | 48.23 | 1247.917 | 518 |
| | Target | 759.3 | 7037.01 | 3436.1 | 58.10 | 1247.917 | 518 |
| ULPGC | Data | 1352.9 | 11720.77 | 2584 | 44.35 | 1360.731 | 425 |
| | Target | 840.1 | 6950.12 | 4149.4 | 50.26 | 1360.731 | 425 |
| USAL | Data | 1917.6 | 16105.84 | 5485 | 49.77 | 1928.015 | 1121 |
| | Target | 1303.9 | 10572.67 | 7060.3 | 53.17 | 1928.015 | 1121 |
| UCAR | Data | 1455.4 | 12564.66 | 6016 | 41.16 | 2138.44 | 465 |
| | Target | 1050.7 | 9582.78 | 6016 | 55.92 | 2138.44 | 754.7 |

Table 4. Actual data and closest targets on $1^{st}$-level efficient frontier

| | | INTERMEDIATE STEP | | | | | | FINAL TARGETS | | | | | |
|---|---|---|---|---|---|---|---|---|---|---|---|---|---|
| | | ASTAFF | EXPEND | ARTICLES | %Q1 | INCOMES | THESIS | ASTAFF | EXPEND | ARTICLES | %Q1 | INCOMES | THESIS |
| | Actual | 1392.5 | 10923.17 | 2546 | 48.23 | 1247.92 | 518 | | | | | | |
| UCA | α=1 | 911.8 | 7926.25 | 2801.0 | 50.22 | 1247.92 | 518 | 759.3 | 7037.01 | 3436.1 | 58.10 | 1247.92 | 518 |
| | α=0.75 | 911.8 | 7926.25 | 2801.0 | 50.22 | 1247.92 | 518 | 759.3 | 7037.01 | 3436.1 | 58.10 | 1247.92 | 518 |
| | α=0.5 | 971.1 | 8873.53 | 3410.6 | 55.96 | 1247.92 | 518 | 759.3 | 7037.01 | 3436.1 | 58.10 | 1247.92 | 518 |
| | α=0.25 | 952.9 | 8769.47 | 4205.9 | 49.45 | 1247.92 | 518 | 831.3 | 6930.70 | 4205.9 | 49.4 | 1247.92 | 518 |
| | α=0 | 1031.1 | 9483.34 | 4782.2 | 57.7 | 1247.92 | 721.0 | 1031.2 | 8325.16 | 4782.2 | 58.77 | 1247.92 | 721.0 |
| | Actual | 1352.9 | 11720.77 | 2584 | 44.35 | 1360.73 | 425 | | | | | | |
| ULPGC | α=1 | 928.2 | 7742.16 | 2747.6 | 50.342 | 1360.73 | 425 | 666.9 | 6883.13 | 3342.6 | 58.27 | 1360.73 | 464.7 |
| | α=0.75 | 927.6 | 7744.92 | 2753.4 | 50.26 | 1360.73 | 425 | 840.1 | 6950.12 | 4149.4 | 50.26 | 1360.73 | 425 |
| | α=0.5 | 912.7 | 7822.26 | 2914.8 | 47.95 | 1360.73 | 425 | 840.1 | 6950.12 | 4149.4 | 50.26 | 1360.73 | 425 |
| | α=0.25 | 1015.1 | 8859.63 | 4213.1 | 50.24 | 1360.73 | 448.2 | 848.4 | 7020.33 | 4213.1 | 50.24 | 1360.73 | 448.2 |
| | α=0 | 1089.7 | 10124.21 | 5171.5 | 57.77 | 1360.73 | 779.0 | 1089.7 | 8793.90 | 5171.5 | 59.16 | 1360.73 | 779.0 |
| | Actual | 1917.6 | 16105.84 | 5485 | 49.77 | 1928.01 | 1121 | | | | | | |
| USAL | α=1 | 1728.4 | 15445.17 | 6694.1 | 49.77 | 2765.53 | 1121 | 1128.1 | 12143.80 | 7231.7 | 59.91 | 2765.53 | 1121 |
| | α=0.75 | 1601.7 | 14739.04 | 6903.3 | 55.32 | 2267.06 | 1121 | 1431.0 | 11477.68 | 7657.9 | 55.32 | 2267.06 | 1121 |
| | α=0.5 | 1570.4 | 15002.47 | 7449.6 | 55.03 | 2221.24 | 1121 | 1413.8 | 11355.38 | 7577.1 | 55.03 | 2221.24 | 1121 |
| | α=0.25 | 1569.3 | 15140.91 | 7614.6 | 54.53 | 2242.49 | 1121 | 1421.8 | 11412.10 | 7614.6 | 55.16 | 2242.49 | 1121 |
| | α=0 | 1569.3 | 15140.91 | 7614.6 | 54.53 | 2242.49 | 1121 | 1421.8 | 11412.10 | 7614.6 | 55.16 | 2242.49 | 1121 |
| | Actual | 1455.4 | 12564.66 | 6016 | 41.16 | 2138.44 | 465 | | | | | | |
| UCAR | α=1 | 1370.1 | 12011.84 | 6016 | 48.03 | 2138.44 | 718.1 | 1050.7 | 9582.78 | 6016 | 55.92 | 2138.44 | 754.7 |
| | α=0.75 | 1370.1 | 12011.84 | 6016 | 48.03 | 2138.44 | 718.1 | 1050.7 | 9582.78 | 6016 | 55.92 | 2138.44 | 754.7 |
| | α=0.5 | 1323.9 | 11974.27 | 6016 | 47.22 | 2138.44 | 720.9 | 1050.7 | 9582.78 | 6016 | 55.92 | 2138.44 | 754.7 |
| | α=0.25 | 1362.0 | 12564.66 | 6016 | 51.85 | 2138.44 | 884.7 | 1286.9 | 10394.88 | 6760.4 | 54.79 | 2138.44 | 884.7 |
| | α=0 | 1362.0 | 12564.66 | 6016 | 51.85 | 2138.44 | 884.7 | 1286.9 | 10394.88 | 6760.4 | 54.79 | 2138.44 | 884.7 |

Table 5. Two-step benchmarking

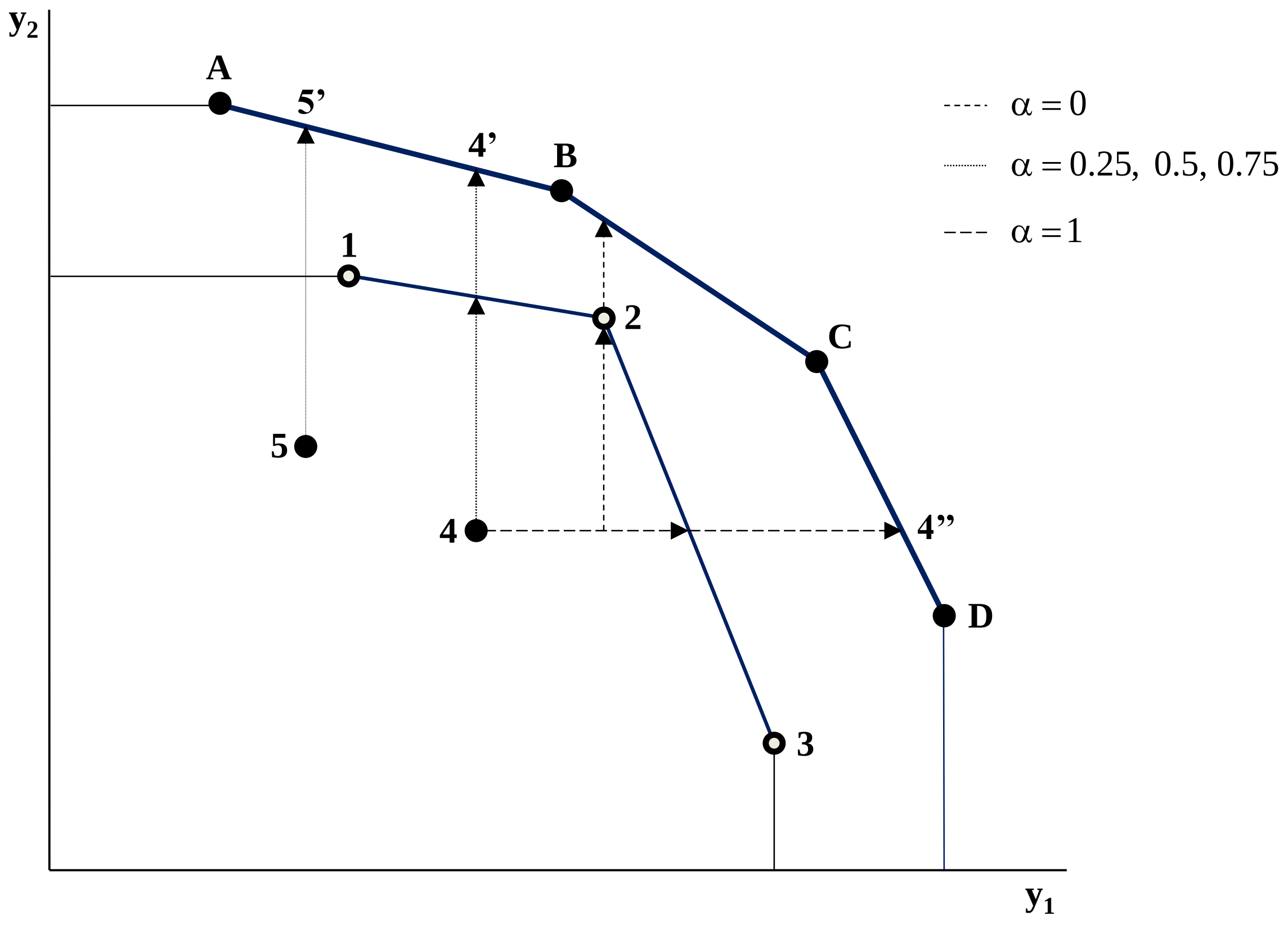

Figure 1. Numerical example